\documentclass{amsart}

\usepackage{hyperref}

\title[Weyl metrisability of projective structures]{Weyl metrisability of two-dimensional\\ projective structures}

\date{May 22, 2013}


\author[T. Mettler]{Thomas Mettler}

\thanks{Research for this article was partially supported by the Swiss National Science Foundation via the grants 200020-121506, 200020-107652, 200021-116165, 200020-124668 and the Postdoctoral Fellowship PBFRP2-133545.}

\address{Department of Applied Mathematics and Theoretical Physics, University of Cambridge, Cambridge, UK.}
\email{t.mettler@damtp.cam.ac.uk}

\newcommand{\J}{\mathfrak{J}}
\renewcommand{\Re}{\operatorname{\mathrm{Re}}}
\renewcommand{\Im}{\operatorname{\mathrm{Im}}}
\renewcommand{\i}{\mathrm{i}}
\renewcommand{\d}{\mathrm{d}}
\newcommand{\D}{\mathrm{D}}
\newcommand{\R}{\mathbb{R}}
\newcommand{\C}{\mathbb{C}}

\newcommand{\cip}{\mathbb{CP}^2\setminus\mathbb{RP}^2}

\newtheorem{theorem}{Theorem}
\newtheorem{lemma}{Lemma}
\newtheorem{corollary}{Corollary}
\newtheorem{proposition}{Proposition}
\theoremstyle{remark}
\newtheorem{remark}{Remark}
\numberwithin{equation}{section}

\begin{document}

\maketitle

\begin{abstract}
We show that on a surface locally every affine torsion-free connection is projectively equivalent to a Weyl connection. First, this is done using exterior differential system theory. Second, this is done by showing that the solutions of the relevant \textsc{PDE} are in one-to-one correspondence with the sections of the `twistor' bundle of conformal inner products having holomorphic image. The second solution allows to use standard results in algebraic geometry to show that the Weyl connections on the two-sphere whose geodesics are the great circles are in one-to-one correspondence with the smooth quadrics without real points in the complex projective plane. 
\end{abstract}


\section{Introduction}

In~\cite{48.0842.02}, Eisenhart and Veblen solve the \textit{Riemannian metrisability problem} for a manifold equipped with a real analytic affine torsion-free connection $\nabla$; i.e.~they determine the necessary and sufficient conditions for $\nabla$ to locally be a Levi-Civita connection or equivalently, the holonomy of $\nabla$ being a subgroup of the orthogonal group. One can also ask to determine the necessary and sufficient conditions for $\nabla$ to be projectively equivalent to a Levi-Civita connection. Recall that two affine connections on a manifold are said to be \textit{projectively equivalent} if they have the same unparametrised geodesics. A projective equivalence class of affine torsion-free connections is called a \textit{projective structure} and will be denoted by $[\nabla]$. Although known since Roger Liouville's initial paper~\cite{21.0317.01} which dates back to 1889, the projective local Riemannian metrisability problem has been solved only recently for real analytic projective structures on surfaces by Bryant, Dunajski and Eastwood~\cite{MR2581355}. A global characterisation of compact Zoll projective surfaces admitting a compatible Levi-Civita connection was given in~\cite{MR1979367}. In~\cite{MR2384718}, the general case is shown to give rise to a linear \textsc{pde} system of finite type. An algorithmic procedure for checking if a given projective structure on a manifold contains a Levi-Civita connection is given in~\cite{arXiv:1003.1469v1} (see also~\cite{arXiv:1101.2069v1}). In~\cite{MR2609304}, it was shown that locally the Riemannian metrisability problem for projective surfaces is equivalent to finding a K\"ahler metric on an associated conformal $4$-manifold of neutral signature. 

There are two problems related to the projective Riemannian metrisability problem that are motivated by two different viewpoints: 

First, the projective Riemannian metrisability problem may be thought of as an inverse problem in the calculus of variations, where one looks for a length functional whose (unparametrised) geodesics are prescribed. However, the functional is constrained to be the length functional of a Riemannian metric. Naturally one might look for a general length functional, more precisely a Finsler metric, whose geodesics are prescribed. This problem is studied in~\cite{alvarezberckfinslermetri} and it is shown that locally on a surface every projective structure (or more generally path geometry) is Finsler metrisable. 

Second, the projective Riemann metrisability problem may be thought of as looking for a connection $\nabla$ in a projective equivalence class $[\nabla]$, whose parallel transport maps are linear \textit{isometries} for some Riemannian metric $g$. From this viewpoint one might also ask for existence of a connection $\nabla \in [\nabla]$, whose parallel transport maps are merely linear \textit{conformal maps} for some conformal structure $[g]$.  It is this latter problem we investigate in this article. More precisely, we study the (projective) \textit{Weyl metrisability problem}, i.e.~the problem of finding an affine torsion-free connection preserving a conformal structure, a so-called \textit{Weyl connection}, whose unparametrised geodesics are prescribed by some projective structure $[\nabla]$. 

Weyl connections were introduced by Weyl~\cite{48.0844.04} as an attempt to unify gravity and electromagnetism and are nowadays mainly studied in the context of the Einstein-Weyl equations in dimensions $d\geq 3$ (see~\cite{MR1807082,MR1199072} and references therein). Also, in~\cite{MR1801870} Wojtkowski observed a relation between Weyl connections and isokinetic dynamics as introduced by Hoover~\cite{hoover} and discussed by Gallavotti and Ruelle in the context of non equilibrium statistical mechanics~\cite{MR1489572}. 

This article is organised as follows. In \S2 we use Cartan's projective connection~\cite{MR1504846} and the theory of exterior differential systems~\cite{MR1083148,MR2003610} to show that locally every smooth projective structure $[\nabla]$ on a surface is Weyl metrisable. In \S3 we characterise the complex structure on the total space of the `twistor bundle'~\cite{MR728412,MR812312} of conformal inner products $\mathcal{C}(M)\to M$ over an oriented surface $M$ in terms of Cartan's projective connection. We use this characterisation to prove the main result: A conformal structure $[g]$ on $M$ is preserved by a $[\nabla]$-representative if and only if~$[g] : M \to \mathcal{C}(M)$ has holomorphic image. As an application of the main result we show in \S4 that the Weyl connections on the $2$-sphere whose geodesics are the great circles are in one-to-one correspondence with the smooth quadrics $C \subset \mathbb{CP}^2$ without real points.

The reader should note that in a certain sense the main results of this article generalise to higher dimensions in the context of Segre structures, see~\cite{MR3043749} for further details.      

\section{An EDS solution for local Weyl metrisability} 

In this section we will use the theory of exterior differential systems (\textsc{eds}) to show that locally every smooth affine torsion-free connection on a surface is projectively equivalent to a Weyl connection. The notation and terminology for \textsc{eds} are chosen to be consistent with~\cite{MR1083148}.  

\subsection{Cartan's projective connection} Weyl showed~\cite{48.0844.04} that two affine connections with the same torsion $\bar{\nabla}$ and $\nabla$ on a smooth manifold $N$ are projectively equivalent if and only if there exists a (unique) $1$-form $\varepsilon \in \mathcal{A}^1(N)$ such that 
\begin{equation}\label{proequi}
\bar{\nabla}_{X}Y-\nabla_{X}Y=\varepsilon(X)\,Y+\varepsilon(Y)\,X
\end{equation}
for every pair of vector fields $X,Y$ on $N$. In more geometric terms, \eqref{proequi} means that the parallel transports of projectively equivalent connections along any curve agree, when thought of as maps between projective space, thus justifying the name projective structure.

As an application of his method of equivalence, Cartan has shown how to associate a parabolic Cartan geometry to a manifold equipped with a projective structure $[\nabla]$. We will only state Cartan's result for oriented projective surfaces, i.e.~two-dimensional, connected, oriented, $C^{\infty}$-manifolds equipped with a $C^{\infty}$ projective structure. For the general case the reader can consult Cartan's original paper~\cite{MR1504846} or~\cite{MR0159284} for a more modern exposition (see~\cite{MR2532439} for background on parabolic Cartan geometries). 

Let $H \subset \mathrm{SL}(3,\mathbb{R})$ be the Lie group of matrices of the form 
$$
H=\left\{\left(\begin{array}{cc} \det(a)^{-1} & b \\ 0 & a\end{array}\right)\;\Bigg\vert\; a\in \mathrm{GL}^+(2,\mathbb{R}), \; b^t \in \mathbb{R}^2\right\}. 
$$
The elements of $H$ will be denoted by $h_{a,b}$. 
\begin{theorem}[Cartan,~\cite{MR1504846}] Let $(M,[\nabla])$ be an oriented projective surface. Then there exists a Cartan geometry $(\pi : B \to M,\theta)$ of type $(\mathrm{SL}(3,\R),H)$ consisting of a right principal $H$-bundle $\pi : B \to M$ and a Cartan connection $\theta \in \mathcal{A}^1(B,\mathfrak{sl}(3,\R))$  with the following properties:  
\begin{itemize}
\item[(i)] Writing $\theta=(\theta^i_j)_{i,j=0,1,2}$,
the leaves of the foliation defined by $\left\{\theta^2_0,\theta^2_1\right\}^{\perp}$ project to geodesics on $M$ and the $\pi$-pullback of every positive volume form on $M$ is a positive multiple of $\theta^1_0\wedge\theta^2_0$.  
\item[(ii)] The \textit{curvature} $2$\textit{-form} $\Theta=\d\theta+\theta\wedge\theta$ satisfies 
\begin{equation}\label{structheta}
\Theta=\left(\begin{array}{ccc} 0 & L_1\;\theta^1_0\wedge\theta^2_0 & L_2\;\theta^1_0\wedge\theta^2_0 \\ 0&0&0 \\ 0&0&0\end{array}\right)
\end{equation}
for some smooth functions $L_i: B \to \mathbb{R}$. 
\end{itemize}
\end{theorem}
Recall that a $2$-\textit{frame at} $p \in M$ is a $2$-jet $j^2_0\varphi$ of a local diffeomorphism $\varphi : U_{0} \to M$ which is defined in a neighbourhood of $0 \in \mathbb{R}^2$ and satisfies $\varphi(0)=p$. The fibre of $\pi : B \to M$ at $p \in M$ precisely consists of those $2$-frames $j^2_0\varphi$ at $p$ for which $\varphi$ is orientation preserving at $0$ and for which $\varphi^{-1}\circ \gamma$ has vanishing curvature at $0$ for every $[\nabla]$-geodesic $\gamma$ through $p$. The Lie group $\tilde{H}$ of $2$-jets of orientation preserving linear fractional transformations $f_{a,b}$ 
$$
x \mapsto \frac{a \cdot x}{1+b\cdot x}, \quad b^t \in \mathbb{R}^2, \;a \in \mathrm{GL}^+(2,\mathbb{R})
$$ 
acts smoothly from the right on $B$ by $j^2_0\varphi \cdot j^2_0f_{a,b}=j^2_0\left(\varphi \circ f_{a,b}\right)$. Note that $H$ and $\tilde{H}$ are isomorphic via the map $h_{a,b} \mapsto j^2_0f_{\tilde{a},\tilde{b}}$ where $\tilde{a}=\det(a)a$ and $\tilde{b}=\det(a)b$. Henceforth we will use this identification whenever needed. 

\subsection{Coordinate sections of Cartan's bundle} Let $x=(x^1,x^2) : U \to \R^2$ be local orientation preserving coordinates on $M$ and $\Gamma^i_{kl} : U \to \R$ denote the Christoffel symbols of a representative of $[\nabla]$ with respect to the coordinates $x$. Then the functions 
$$
\Pi^i_{kl}=\Gamma^i_{kl}-\frac{1}{3}\left(\delta^i_k \sum_{j}\Gamma^j_{jl}+\delta^i_l\sum_{j}\Gamma^j_{jk}\right)
$$
are projective invariants in the sense that they do not depend on the representative chosen to compute them, but only on $x$. Locally $[\nabla]$ can be recovered from the projective invariants $\Pi^i_{kl}$ by defining them to be the Christoffel symbols with respect to $x$ of an affine torsion-free connection $\nabla^{\prime}$, which is well defined on $U$ and projectively equivalent to $[\nabla]$. Consequently, two affine torsion-free connections on $M$ are projectively equivalent, if and only if their Christoffel symbols give rise to the same functions $\Pi^i_{kl}$. Associated to the coordinates $x$ is a \textit{coordinate section} $\sigma_x : U \to B$ which assigns to every point $p \in U$ the $2$-frame $j^2_0\varphi \in B$ at $p$ defined by 
\begin{equation}\label{2framecartan}
\varphi(0)=p, \quad \partial_k(x\circ \varphi)^i(0)=\delta^i_k, \quad \partial_k\partial_l(x\circ \varphi)^i(0)=-\Pi^i_{kl}(p). 
\end{equation}
This section indeed does take values in $B$ as can be shown with a simple computation. The section $\sigma_x: U \to B$ satisfies
\begin{equation}\label{coordinatesection}
\sigma_x^*\theta^0_0=0, \quad \sigma_x^*\theta^1_0=\d x^1, \quad \sigma_x^*\theta^2_0=\d x^2,  
\end{equation}
thus the structure equations \eqref{structheta} yield
$$
\aligned
0&=\sigma_x^*\d\theta^0_0=-\sigma_x^*\theta^0_1\wedge \d x^1-\sigma_x^*\theta^0_2\wedge \d x^2,\\
0&=\sigma_x^*\d\theta^1_0=-\sigma_x^*\theta^1_1\wedge \d x^1-\sigma_x^*\theta^1_2\wedge \d x^2,\\
0&=\sigma_x^*\d\theta^2_0=-\sigma_x^*\theta^2_1\wedge \d x^1-\sigma_x^*\theta^2_2\wedge \d x^2.\\
\endaligned
$$
Since $\sigma_x^*(\theta^0_0)=-\sigma_x^*(\theta^1_1+\theta^2_2)=0$ holds, Cartan's lemma implies that there exist functions $\kappa_0,\kappa_1,\kappa_2,\kappa_3$ and $\zeta_1,\zeta_2,\zeta_3$ on $U$ such that $\sigma_x^*\theta=\eta_x$ where
\begin{equation}\label{pbcs}
\eta_x=\left(\begin{array}{crr} 0& \zeta_1 \d x^1+\zeta_2 \d x^2 & \zeta_2 \d x^1+\zeta_3 \d x^2\\ \d x^1 & -\kappa_1 \d x^1-\kappa_2 \d x^2 & -\kappa_2 \d x^1-\kappa_3 \d x^2  \\ \d x^2 & \kappa_0 \d x^1+\kappa_1 \d x^2 & \kappa_1 \d x^1+\kappa_2 \d x^2\end{array}\right).
\end{equation}
The $\theta$-structure equations then imply that the functions $\zeta_i$ and $\kappa_j$ satisfy the relations
\begin{equation}\label{relzetakappa}
\aligned
\zeta_1&=\frac{\partial \kappa_1}{\partial x^1}-\frac{\partial \kappa_0}{\partial x^2}+2\kappa_1^2-2\kappa_0\kappa_2,\\
\zeta_2&=\frac{\partial \kappa_2}{\partial x^1}-\frac{\partial \kappa_1}{\partial x^2}+\kappa_1\kappa_2-\kappa_0\kappa_3,\\
\zeta_3&=\frac{\partial \kappa_3}{\partial x^1}-\frac{\partial \kappa_2}{\partial x^2}+2\kappa_2^2-2\kappa_1\kappa_3.\\
\endaligned
\end{equation}
In terms of the $\Pi^i_{kl}$, the functions $\kappa_i : U \to \R$ can be expressed as
$$
\kappa_0=\Pi^2_{11}, \quad \kappa_1=\Pi^2_{12},\quad \kappa_2=\Pi^2_{22}, \quad \kappa_3=-\Pi^1_{22}.
$$
The \textit{coefficients} $\Pi^i_{kl}$ of Cartan's projective connection were discovered independently of Cartan's work from the invariant theoretic view point by Thomas~\cite{51.0569.03}. They generalise to the $n$-dimensional case by replacing $3$ with $n+1$. 

Note that the coordinate section $\sigma_x : U \to B$ is an integral manifold of the \textsc{eds} with independence condition $(\mathcal{I},\Omega)$ where 
$$
\mathcal{I}=\left\langle \theta^0_0,\d\theta^1_0,\d\theta^2_0\right\rangle
$$
and $\Omega=\theta^1_0\wedge\theta^2_0$. 
Conversely, if $s : U \to B$ is an integral manifold of the \textsc{eds} $(\mathcal{I},\Omega)$, then $\pi \circ s : U \to M$ is a local diffeomorphism, thus locally, without losing generality, we can assume that $U \subset M$ and $s$ is a section of $\pi : B \to M$. Assume $U$ is simply connected, then
$$
s^*\d\theta^1_0=s^*\d\theta^2_0=0 
$$ 
implies that there exist functions $(x^1,x^2) : U \to \R^2$ such that
$$
s^*\theta^1_0=\d x^1, \quad s^*\theta^2_0=\d x^2.
$$
Now $s^*(\theta^1_0\wedge\theta^2_0)=\d x^1\wedge \d x^2 > 0$ implies that $x=(x^1,x^2) : U \to \R^2$ is a local orientation preserving coordinate system on $M$ satisfying $\sigma_x=s$.  
 
\subsection{Coordinate sections and Weyl metrisability} 
Note that an affine torsion-free connection $\nabla$ preserves a conformal structure $[g]$, if and only if for some (and hence any) Riemannian metric $g \in [g]$, there exists a $1$-form $\beta$, so that 
\begin{equation}\label{hans}
\nabla g=2 \beta \otimes g. 
\end{equation}
The affine torsion-free connections preserving $[g]$ are called \textit{Weyl connections} for $[g]$. Given a Riemannian metric $g$ and a $1$-form $\beta$, the affine torsion-free connection $\D^{g,\beta}$ defined by
$$
(X,Y) \mapsto \D^g_{X}Y+g(X,Y)\beta^{\sharp}-\beta(X)Y-\beta(Y)X,
$$
is the unique Weyl connection for $[g]$ solving \eqref{hans}. Here $\D^g$ denotes the Levi-Civita connection of $g$ and $\beta^{\sharp}$ the $g$-dual vector field to $\beta$.  

On the total space $B$ of the Cartan geometry of an oriented projective surface $(M,[\nabla])$ consider the \textsc{eds} with independence condition $(\mathcal{I},\Omega)$ defined by
\begin{equation}\label{firstedsweylmetri}
\mathcal{I}=\left\langle \theta^0_0, \d\theta^1_0, \d\theta^2_0,\theta^1_0\wedge(3\,\theta^2_1+\theta^1_2),\theta^2_0\wedge(\theta^2_1+3\,\theta^1_2) \right \rangle, \quad \Omega =\theta^1_0\wedge\theta^2_0.
\end{equation} 
This \textsc{eds} is of interest due to the following: 
\begin{lemma}\label{edsweylmetri}
Let $\nabla$ be a Weyl connection on the oriented surface $M$ and $(\pi : B\to M,\theta)$ the Cartan geometry associated to $[\nabla]$. Then in a neighbourhood $U_p$ of every point $p \in M$, there exists a coordinate section $\sigma_x : U_p \to B$ which is an integral manifold of $(\mathcal{I},\Omega)$. Conversely, let $(M,[\nabla])$ be an oriented projective surface with Cartan geometry $(\pi : B \to M,\theta)$. Then every coordinate section $\sigma_x : U\subset M \to B$ which is an integral manifold of $(\mathcal{I},\Omega)$ gives rise to a Weyl connection on $U$ which is projectively equivalent to $\nabla$.      
\end{lemma}

\begin{proof}
Let $\nabla$ be a Weyl connection for $[g]$ and $(\pi : B \to M,\theta)$ the Cartan geometry associated to $[\nabla]$. For a given point $p \in M$ let $x=(x^1,x^2) : U_p \to \mathbb{R}^2$ be local $p$-centred coordinates which are orientation preserving and isothermal for $[g]$. Then there exist smooth functions $r_i: U \to \R$ such that
$$
\nabla (x^*g_E)=2(r_1\d x^1+r_2 \d x^2)\otimes x^*g_E,
$$
where $g_E$ denotes the Euclidean standard metric on $\R^2$. Now a simple computation shows that the projective invariants $\kappa_i : U \to \R$ of the projective structure $[\nabla]$, defined with respect to $x$, satisfy the relations 
\begin{equation}\label{weylcond}
\kappa_0=3\kappa_2=r_2, \quad 3\kappa_1=\kappa_3=-r_1. 
\end{equation}
It follows with \eqref{pbcs} that the associated coordinate section $\sigma_x$ satisfies 
$$
\aligned
\sigma_x^*\left(\theta^1_0\wedge(3\,\theta^2_1+\theta^1_2)\right)&=\left(3\kappa_1-\kappa_3\right)\d x^1\wedge \d x^2=0,\\
\sigma_x^*\left(\theta^2_0\wedge(\theta^2_1+3\,\theta^1_2)\right)&=\left(3\kappa_2-\kappa_0\right)\d x^1\wedge \d x^2=0,\\
\endaligned
$$
thus showing that $\sigma_x$ is an integral manifold of $(\mathcal{I},\Omega)$. Conversely, let $(M,[\nabla])$ be an oriented projective surface with Cartan geometry $(\pi : B \to M,\theta)$. Suppose $\sigma_x : U \to B$ is a coordinate section and an integral manifold of $(\mathcal{I},\Omega)$. Then the projective invariants $\kappa_i : U\to \R$ with respect to $x$ satisfy \eqref{weylcond} and thus the Weyl connection $\D^{g,\beta}$ defined on $U$ by the pair  
$$
g=x^*g_E,\quad \beta=-\kappa_3\d x^1+\kappa_0\d x^2,
$$
is projectively equivalent to $\nabla$. 
\end{proof}

Lemma \ref{edsweylmetri} translates the Weyl metrisability problem into finding integral manifolds of the \textsc{eds} $(\mathcal{I},\Omega)$.  For the application of the theory of exterior differential systems it is more convenient to work with a linear Pfaffian system. Let $A=B\times \mathbb{R}^2$ and denote by $a_i : A \to \R$ the projection onto the $i$-th coordinate of $\R^2$ and by $\tau : A \to B$ the canonical projection. On $A$ define 
$$
\omega^1=\tau^*\theta^1_0, \quad \omega^2=\tau^*\theta^2_0,
$$   
and
$$
\aligned
\vartheta^1&=\tau^*\theta^1_1+a_1\omega^1+a_2\omega^2,\\
\vartheta^2&=\tau^*\theta^2_2-a_1\omega^1-a_2\omega^2,\\
\vartheta^3&=\tau^*\theta^1_2+a_2\omega^1+3a_1\omega^2,\\
\vartheta^4&=\tau^*\theta^2_1-3a_2\omega^1-a_1\omega^2.
\endaligned
$$
Then it is straightforward to check that the integral manifolds of the \textsc{eds}
$$
\mathcal{I}^{\prime}=\left\langle \vartheta^1,\vartheta^2,\vartheta^3,\vartheta^4\right\rangle, \quad \Omega^{\prime}=\omega^1\wedge\omega^2,
$$
are in one-to-one correspondence with the integral manifolds of $(\mathcal{I},\Omega)$. We follow the strategy explained in~\cite[Chapter IV, \S2]{MR1083148} to prove: 
\begin{theorem}\label{edsweylmetrimain}
The \textsc{eds} $(\mathcal{I}^{\prime},\Omega^{\prime})$ is locally equivalent to a determined first order elliptic \textsc{pde} system for $4$ real-valued functions of $2$ variables. In particular locally every projective surface is Weyl metrisable.  
\end{theorem}
\begin{proof}
Let
$$
\aligned
I&=\text{span}\left\{\vartheta^1,\vartheta^2,\vartheta^3,\vartheta^4\right\} \subset T^*A,\\
J&=\text{span}\left\{\vartheta^1,\vartheta^2,\vartheta^3,\vartheta^4,\omega^1,\omega^2\right\} \subset T^*A,\\
L&=J/I\simeq\text{span}\left\{\omega^1,\omega^2\right\} \subset T^*A,\\
\endaligned
$$
where span means linear combinations with coefficients in $C^{\infty}(A,\R)$. It is easily verified that $J$ is a Frobenius system, in particular $(\mathcal{I}^{\prime},\Omega^{\prime})$ is locally equivalent to a first order \textsc{pde} system
\begin{equation}\label{pdesys}
F^b(y^i,z^a,\partial z^a/\partial y^i)=0
\end{equation}
for $4=\text{rank}\,I$ real-valued functions $z^a$ of $2=\text{rank}\, L$ variables $y^i$~\cite[Prop. 5.10, Chapter IV]{MR1083148}. Let $\lambda : P=\mathbb{P}(L) \to A$ denote the  bundle which is obtained by projectivisation of the (quotient) vector-bundle $L \to A$.  We have the structure equations 
$$
\d\vartheta^i=\sum_{k=1}^2\varphi^i_k\wedge \omega^k, \quad \text{mod}\; \vartheta_l,\; l=1,\ldots,4
$$
with
\begin{equation}\label{symbchar}
\varphi^i_k=\left(\begin{array}{cc}\varphi_1 & \varphi_2 \\ \varphi_3 & \varphi_4 \\ -3\varphi_2 & 2\varphi_3+\varphi_1 \\ 2\varphi_2+\varphi_4 & -3\varphi_3\end{array}\right)
\end{equation}
for $4$ linearly independent $1$-forms $\varphi_l$ on $A$. Therefore we have $4$ non-trivial (symbol) relations for the $8$ entries of $\varphi^i_k$ and since $s_0=\text{rank}\,I=4$, the linear Pfaffian system $(\mathcal{I}^{\prime},\Omega^{\prime})$ and the corresponding \textsc{pde} system \eqref{pdesys} are determined. Moreover straightforward computations using \eqref{symbchar} show that the characteristic variety $\Xi \subset P$ of $(\mathcal{I}^{\prime},\Omega^{\prime})$ at $a \in A$ is given by
$$
\Xi_a=\left\{\left[\xi_1\omega^1(a)+\xi_2\omega^2(a)\right] \in \lambda^{-1}(a) \, \Big| \left((\xi_1)^2+(\xi_2)^2\right)^2=0 \, \right\},
$$ 
thus $\Xi$ is empty. This shows that $(\mathcal{I}^{\prime},\Omega^{\prime})$ and the corresponding \textsc{pde} system \eqref{pdesys} are elliptic. It follows with standard results in elliptic \textsc{pde} theory (see~\cite[p.~15]{MR0450755}) that \eqref{pdesys} has smooth local solutions. Using Lemma \ref{edsweylmetri} we conclude that locally every smooth projective surface is Weyl metrisable.  
\end{proof}

\section{A complex geometry solution for Weyl metrisability}

In this section we will show that the Weyl metrisability problem for an oriented projective surface is globally equivalent to finding a section of the bundle of conformal inner products with holomorphic image. 

\subsection{The bundle of conformal inner products} Recall that a \textit{conformal inner product} on a real vector space $V$ is an equivalence class $[b]$ of inner products on $V$, where two inner products are called equivalent if one is a positive multiple of the other. 

Let $N$ be a manifold of even dimension $2n$ and let $\mathcal{F}(N) \to N$ be the right principal $\mathrm{GL}(2n,\R)$-bundle of $1$-frames over $N$. We embed $\mathrm{GL}(n,\C)$ as a closed subgroup of $\mathrm{GL}(2n,\R)$. Let $\mathcal{F}(N)/\mathrm{GL}(n,\C) \to N$ be the  bundle whose fibre at $p\in N$ consists of the complex structures on $T_pN$. It was observed in~\cite{MR728412,MR812312} that the choice of an affine connection $\nabla$ on $N$ induces an almost complex structure $\J$ on $\mathcal{F}(N)/\mathrm{GL}(n,\C)$. If $\nabla$ is torsion-free, then $\J$ is integrable if and only if the Weyl projective curvature tensor of $\nabla$ vanishes. In fact, $\J$ only depends on the projective equivalence class of $\nabla$. In the case where $N$ is oriented, this almost complex structure $\J$ restricts to become an almost complex structure on the subbundle $\mathcal{F}^+(N)/\mathrm{GL}(n,\C)$ where $\mathcal{F}^+(N) \to N$ denotes bundle of positively oriented frames.

For the case of an oriented surface $M$, the fibre of the bundle $\rho : \mathcal{C}(M)=\mathcal{F}^+(M)/\mathrm{GL}(1,\C)\to M$ at $p\in M$ may be identified with the space of conformal inner products on $T_pM$. Consequently, a conformal structure on $M$ may also be thought of as a section of the \textit{bundle of conformal inner products} $\rho : \mathcal{C}(M) \to M$. Note that in two dimensions the Weyl projective curvature tensors vanishes identically for every projective structure $[\nabla]$. It follows that the almost complex structure $\J$ is always integrable.  

\subsection{The complex surface $\mathcal{C}(M)$ and Cartan's connection} In this subsection we will characterise the complex structure on $\mathcal{C}(M)$ in terms of the Cartan geometry $(\pi : B \to M, \theta)$ associated to $[\nabla]$. To this end let $C\subset H$ be the closed Lie subgroup consisting of elements $h_{a,b}$ with $a \in \mathrm{GL}(1,\C)$ where we identify $\mathrm{GL}(1,\C)$ with the non-zero $2$-by-$2$ matrices of the form 
$$
\left(\begin{array}{rr} x & -y \\ y & x \end{array}\right).
$$ Consider the smooth map 
$$
\nu : B \to \mathcal{C}(M), \; j^2_0\varphi \mapsto \left[(\varphi_*g_E)_{\varphi(0)}\right].
$$ 
\begin{lemma}
The map $\nu : B \to \mathcal{C}(M)$ makes $B$ into a right principal $C$-bundle over $\mathcal{C}(M)$.  
\end{lemma}

\begin{proof}
Since $C$ is a closed Lie subgroup of $H$, it is sufficient to show that $\nu$ is a smooth surjection whose fibres are the $C$-orbits. Clearly $\nu$ is smooth and surjective. Suppose $\nu(j^2_0\varphi)=\nu(j^2_0\tilde{\varphi})$ for some elements $j^2_0\varphi, j^2_0\tilde{\varphi}  \in B$. Then these two elements are in the same fibre of $\pi : B \to M$, hence there exists $h_{a,b} \in H$ such that $j^2_0\tilde{\varphi}=j^2_0\varphi \cdot h_{a,b}$ and 
\begin{equation}\label{gliich}
c \left(\varphi_*g_E\right)_{\varphi(0)}=\left(\left(\varphi \circ f_{a,b} \right)_*g_E\right)_{\varphi(0)} 
\end{equation}
which is equivalent to
$$
c(g_E)_0=(\tilde{a}_*g_E)_0,
$$
where $\tilde{a} \in \mathrm{GL}^+(2,\R)$ is the linear map $x \mapsto (\det a)\, a\cdot x$ and $c\in \mathbb{R}^+$. This is equivalent to $\tilde{a}$ being in $\mathrm{GL}(1,\C)$ or $h_{a,b} \in C$. In other words, the $\nu$ fibres are the $C$-orbits. 
\end{proof}

We will now use the forms $\alpha_1=\theta^1_0+\i\theta^2_0$ and $\alpha_2=(\theta^1_2+\theta^2_1)+\i(\theta^2_2-\theta^1_1)$ to define an almost complex structure on $\mathcal{C}(M)$. Note that the forms $\alpha_1, \alpha_2$ are $\nu$-semibasic, i.e.~$\alpha_i(X)$ vanishes for every vector field $X\in\mathfrak{X}(B)$ which is tangent to the fibres of $\nu$. 

\begin{proposition}\label{twistorcartan}
There exists a unique complex structure $J$ on $\mathcal{C}(M)$ such that a complex valued $1$-form $\mu \in \mathcal{A}^1(\mathcal{C}(M),\mathbb{C})$ is of type $(1,\! 0)$ if and only if $\nu^*\mu$ is a linear combination of $\left\{\alpha_1,\alpha_2\right\}$ with coefficients in $C^{\infty}(B,\mathbb{C})$. 
\end{proposition}
\begin{proof}
Let $T^i_j$ denote the vector fields dual to the coframing $\theta^i_j$. For $\xi \in T\mathcal{C}(M)$ define 
\begin{multline*}
J(\xi)=\nu^{\prime}\left(-\theta^2_0(\tilde{\xi})T^1_0+\theta^1_0(\tilde{\xi})T^2_0-\frac{1}{2}(\theta^2_2-\theta^1_1)(\tilde{\xi})(T^1_2+T^2_1)+\right.\\\left.+\frac{1}{2}(\theta^1_2+\theta^2_1)(\tilde{\xi})(T^2_2-T^1_1)\right),
\end{multline*}
where $\tilde{\xi} \in TB$ satisfies $\nu^{\prime}(\tilde{\xi})=\xi$. Any other vector in $TB$ which is mapped to $\xi$ under $\nu^{\prime}$ is of the form $(R_c)^{\prime}(\tilde{\xi})+\chi$ for some $c\in C$ and $\chi \in \ker\nu^{\prime}$. Using the identities 
$$(R_c)^*\theta=c^{-1}\,\theta\, c, \quad \nu \circ R_c=\nu, \quad c \in C$$ and the fact that $\alpha_1,\alpha_2$ are $\nu$-semibasic, it follows from straightforward computations that $J$ is a well defined almost complex structure on $\mathcal{C}(M)$ which has all the desired properties. Moreover the structure equations \eqref{structheta} imply 
\begin{equation}
\aligned
\d\alpha_1=&\left(-3\theta^2_2+i\left(2\theta^1_2+\theta^2_1\right)\right)\wedge\alpha_1+\left(-\theta^2_0+2i\theta^1_0\right)\wedge\alpha_2,\\
\d\alpha_2=&\left(\theta^0_2-i\theta^0_1\right)\wedge\alpha_1+i\left(\theta^1_2-\theta^2_1\right)\wedge\alpha_2,\\
\endaligned
\end{equation}
and hence, by Newlander-Nirenberg~\cite{MR0088770}, $J$ is integrable. Clearly such a complex structure is unique.  
\end{proof}
In fact it is not hard to show that every $\rho$-fibre admits the structure of a Riemann surface biholomorphic to the unit disk $D^2$ such that the canonical inclusion into $\mathcal{C}(M)$ is a holomorphic embedding. 

\subsection{The compatibility problem and holomorphic curves} In this subsection we will relate holomorphic curves in $\mathcal{C}(M)$ to the Weyl metrisability problem. We will use the following lemma whose prove is elementary and thus omitted.  
\begin{lemma}\label{sub}
Let $(X,\! J)$ be a complex surface, $\mu_1,\mu_2 \in \mathcal{A}^1(X,\mathbb{C})$ a basis for the $(1,\! 0)$-forms of $J$ and $f : \Sigma \to X$ a $2$-submanifold with 
$$f^*(\Re(\mu_1)\wedge\Im(\mu_1))\neq 0.$$ 
Then $f : \Sigma \to X$ is a holomorphic curve if and only if $f^*(\mu_1\wedge\mu_2)=0$. Moreover through every point $p \in X$ passes such a holomorphic curve.  
\end{lemma}
\begin{remark}
Here a $2$-submanifold $f : \Sigma \to X$ is called a \textit{holomorphic curve} if $(J\circ f^{\prime})(T_p\Sigma)=f^{\prime}(T_p\Sigma)$ for every $p \in \Sigma$.
\end{remark}
Let $[g]$ be a conformal structure on $(M,[\nabla])$ and $x : U \to \R^2$ local orientation preserving coordinates which are isothermal for $[g]$. Then it is easy to check that the coordinate section $\sigma_x : U \to B$ satisfies $\nu \circ \sigma_x= [g]\vert_U$.

We will now relate conformal structures $[g] : M \to  \mathcal{C}(M)$, which are holomorphic curves in $\mathcal{C}(M)$ to the \textsc{eds} with independence condition $(\mathcal{I},\Omega)$ on $B$ given by
$$
\mathcal{I}=\left\langle \theta^0_0, \d\theta^1_0,\d\theta^2_0,\Re(\alpha_1\wedge\alpha_2),\Im(\alpha_1\wedge\alpha_2) \right\rangle, \quad \Omega=\Re(\alpha_1)\wedge\Im(\alpha_1).
$$
Note that we may write
$$
\left(\theta^1_0\wedge(3\theta^2_1+\theta^1_2)\right)+\i\left(\theta^2_0\wedge(\theta^2_1+3\theta^1_2)\right)=\alpha_1\wedge\alpha_2+3\i\bar{\alpha}_1\wedge\theta^0_0+2\i\d\bar{\alpha}_1
$$
where $\bar{\alpha}_1=\theta^1_0-i\theta^2_0$. It follows that the \textsc{eds} $(\mathcal{I},\Omega)$ equals the \textsc{eds} \eqref{firstedsweylmetri}.  
\begin{lemma}\label{eds} Let $[g]$ be a conformal structure on $(M,[\nabla])$ and $x=(x^1,x^2) : U \to \mathbb{R}^2$ local orientation preserving $[g]$-isothermal coordinates. Then the coordinate section $\sigma_x : U \to B$ is an integral manifold of $(\mathcal{I},\Omega)$ if and only if $[g]\vert_U : U \to \mathcal{C}(M)$ is a holomorphic curve. 
\end{lemma}

\begin{proof}
Let $s : \rho^{-1}(U) \to B$ be a local section of the bundle $\nu : B\to \mathcal{C}(M)$ and let $\mu_1=s^*\alpha_1, \mu_2=s^*\alpha_2$ be a local basis for the $(1,\! 0)$-forms on $\rho^{-1}(U)$. Note that such sections exist, since the principal bundle $H \to H/C$ is trivial. Now 
$$
\nu^*\mu_k=(s \circ \nu)^*\alpha_k=R_t^*\alpha_k,
$$
for some smooth function $t : \pi^{-1}(U) \to C$. Write the elements of $C\subset H$ in the form 
$$
c_{z,w}=\left(\begin{array}{ccr} \vert z \vert^{-2} & \Re(w) & \Im(w) \\ 
0 & \Re(z) & -\Im(z) \\ 
0 & \Im(z) & \Re(z) \end{array}\right),  
$$
for some complex numbers $z \neq 0$ and $w$. Since $\theta$ is a Cartan connection, we have $R_h^*\theta=h^{-1}\theta h,$ for every $h \in H$ which yields together with a short computation
\begin{equation}\label{trafocomplex}
R_{c_{z,w}}^*\left(\begin{array}{c} \alpha_1 \\ \alpha_2\end{array}\right)=\left(\begin{array}{cc} \bar{z}/\vert z \vert^4 & 0 \\ \i\bar{w}/\bar{z} & \bar{z}/z \end{array}\right)\left(\begin{array}{c} \alpha_1 \\ \alpha_2\end{array}\right).
\end{equation}
Using \eqref{trafocomplex} it follows 
$$
\nu^*\mu_1=\lambda_1 \alpha_1, \quad \nu^*\mu_2=\lambda_2 \alpha_1 + \lambda_3 \alpha_2,
$$
for some smooth functions $\lambda_k : \pi^{-1}(U) \to \mathbb{C}$ with $\lambda_1\lambda_3 \neq 0$. Then \eqref{pbcs} and the identity $[g]\vert_U=\nu \circ \sigma_x$ yield
$$
\aligned
\sigma_x^*\left(\Re(\alpha_1)\wedge\Im(\alpha_1)\right)=\d x^1\wedge \d x^2&\neq 0,\\\sigma_x^*\d\alpha_1=\d(\d x^1+\i\d x^2)&=0,\\ \sigma_x^*\theta^0_0&=0,\\ 
\endaligned
$$
and
\begin{equation}\label{pullback}
([g]\vert_U)^*\mu_1=(\nu\circ\sigma_x)^*\mu_1=\sigma_x^*(\lambda_1\alpha_1)=(\lambda_1\circ \sigma_x)(\d x^1+\i\d x^2)\neq 0, 
\end{equation}
which shows that $([g]\vert_U)^*(\Re(\mu_1)\wedge\Im(\mu_1))\neq 0$. Therefore according to Lemma \ref{sub}, $[g]\vert_U : U \to \mathcal{C}(M)$ is a holomorphic curve if and only if 
$$
([g]\vert_U)^*(\mu_1\wedge\mu_2)=(\nu\circ\sigma_x)^*(\mu_1\wedge\mu_2)=(\lambda_1\lambda_3\circ \sigma_x)\,\sigma_x^*(\alpha_1\wedge\alpha_2)=0, 
$$
which finishes the proof. 
\end{proof}
The \textsc{eds} $(\mathcal{I},\Omega)$ precisely governs the Weyl metrisability problem for an oriented projective surface. 
\begin{proposition}\label{equiv2}
Let $[g]$ be a conformal structure on $(M,[\nabla])$. Then the following two statements are equivalent:
\begin{itemize}
\item[(i)] There exists a Weyl connection for $[g]$ on $M$ which is projectively equivalent to $\nabla$.
\item[(ii)] The coordinate section $\sigma_x : U \to B$ associated to any local orientation preserving $[g]$-isothermal coordinate chart $x=(x^1,x^2) : U \to \mathbb{R}^2$ is an integral manifold of $(\mathcal{I},\Omega)$.
\end{itemize} 
\end{proposition} 
\begin{proof}
(i) $\Rightarrow$ (ii): This direction is an immediate consequence of Lemma \ref{edsweylmetri}.\\
(ii)$\Rightarrow$ (i): Let $x : U\to\mathbb{R}^2$, be local orientation preserving isothermal coordinates for $[g]$ and $\sigma_{x} : U \to B$ the corresponding coordinate section which is an integral manifold of $(\mathcal{I},\Omega)$. Fix a $[g]-$representative $g$, then $$g\vert_{U}=e^{2f}x^*g_E$$ for some smooth~$f : U \to \mathbb{R}$. Since $\sigma_x$ is an integral manifold, the projective invariants $\kappa_i : U \to \mathbb{R}$ with respect to $x$ satisfy $3\kappa_1=\kappa_3$ and $3\kappa_2=\kappa_0$. On $U$ define the $1$-form 
$$\beta=-\kappa_3 \d x^1+\kappa_0 \d x^2+\d f,$$ 
then the Weyl connection $\D^{g,\beta}$ on $U$ associated to the pair $(g\vert_U,\beta)$ is projectively equivalent to $\nabla$. Let~$\tilde{x} : \tilde{U} \to \mathbb{R}^2$ be another local orientation preserving isothermal coordinate chart for $[g]$ overlapping with $U$. Writing $g\vert_{\tilde{U}}=e^{2\tilde{f}}\tilde{x}^*g_E$ for some smooth $\tilde{f} : \tilde{U} \to \mathbb{R}$ and $\tilde{\kappa}_i : \tilde{U} \to\mathbb{R}$ for the projective invariants of $[\nabla]$ with respect to $\tilde{x}$, then again the Weyl connection $\D^{g,\tilde{\beta}}$ on $\tilde{U}$ associated to the pair $(g\vert_{\tilde{U}},\tilde{\beta})$ with 
$$\tilde{\beta}=-\tilde{\kappa}_3\d\tilde{x}^1+\tilde{\kappa}_0\d\tilde{x}^2+\d\tilde{f}
$$ is projectively equivalent to $\nabla$. On $U\cap\tilde{U}$ we have 
$$
\left[\D^g+g\otimes \beta^{\sharp}\right]=[\nabla]=\left[\D^g+g\otimes\tilde{\beta}^{\sharp}\right],
$$ 
using Weyl's result \eqref{proequi}, this is equivalent to the existence of a $1$-form $\varepsilon$ on $U\cap\tilde{U}$ such that 
$$
\D^g_{X}Y+g(X,Y)\beta^{\#}=\D^g_{X}Y+g(X,Y)\tilde{\beta}^{\#}+\varepsilon(X)Y+\varepsilon(Y)X
$$
for every pair of vector fields $X,Y$ on $U\cap \tilde{U}$. In particular the choice of a basis of $g$-orthonormal vector fields implies $\varepsilon=0$. Thus~$\beta=\tilde{\beta}$ on $U \cap \tilde{U}$ and therefore, using a coordinate cover, $\beta$ extends to a well defined global $1$-form which proves the existence of a smooth Weyl connection on $M$ which is projectively equivalent to $\nabla$.   
\end{proof}
Summarising the results found so far we have the main 
\begin{theorem}\label{main}
A conformal structure $[g]$ on an oriented projective surface $(M,[\nabla])$ is preserved by a $[\nabla]$-representative if and only if $[g] : M\to \mathcal{C}(M)$ is a holomorphic curve. 
\end{theorem}
\begin{proof} 
This follows immediately from Lemma \ref{eds} and Proposition \ref{equiv2}. 
\end{proof}
\begin{remark}
It is easy to check that for a given projective structure $[\nabla]$ on $M$ every holomorphic curve $[g] : M \to  \mathcal{C}(M)$ determines a unique Weyl connection which is projectively equivalent to $\nabla$. Theorem \ref{main} therefore gives a one-to-one correspondence between the Weyl connections on an oriented surface whose unparametrised geodesics are prescribed by a projective structure $[\nabla]$ and sections of the fibre bundle $\rho : \mathcal{C}(M) \to M$ which are holomorphic curves. Note also that a conformal structure $[g]$ on an oriented surface determines a unique complex structure whose holomorphic coordinates are given by orientation preserving isothermal coordinates for $[g]$. It follows that a conformal structure $[g]$ on an oriented projective surface $(M,[\nabla])$ is preserved by a $[\nabla]$-representative if and only if $[g] : M\to \mathcal{C}(M)$ is holomorphic with respect to the complex structure on $\mathcal{C}(M)$ and the complex structure on $M$ induced by $[g]$ and the orientation.
\end{remark}
\begin{corollary}\label{localweylmetri}
An affine torsion-free connection $\nabla$ on a surface $M$ is locally projectively equivalent to a Weyl connection. 
\end{corollary}
\begin{proof}
Since the statement is local we may assume that $M$ is oriented. Let $(\pi : B \to M,\theta)$ be the Cartan geometry associated to $[\nabla]$. For a given point $p \in M$, choose $q \in \mathcal{C}(M)$ with $\rho(q)=p$ and a coordinate neighbourhood $U_p$. Let $\mu_1,\mu_2$ be a basis for the $(1,\! 0)$-forms on $\rho^{-1}(U_p)$ as constructed in Lemma \ref{eds}. Using Lemma \ref{sub} there exists a complex $2$-submanifold $f : \Sigma \to \rho^{-1}(U_p)$ passing through $q$ for which $f^*(\Re(\mu_1)\wedge\Im(\mu_1))\neq 0$. Since the $\pi : B \to M$ pullback of a volume form on $M$ is a nowhere vanishing multiple of $\theta^1_0\wedge\theta^2_0$, the $\rho$ pullback of a volume form on $U_p$ is a nowhere vanishing multiple of $\Re(\mu_1)\wedge\Im(\mu_1)$ and hence $\rho \circ f : \Sigma \to U_p$ is a local diffeomorphism. Composing $f$ with the locally available inverse of this local diffeomorphism one gets a local section of the bundle of conformal inner products which is defined in a neighbourhood of $p$ and which is a holomorphic curve. Using Theorem \ref{main} it follows that $\nabla$ is locally projectively equivalent to a Weyl connection.\end{proof}

\section{The flat case}

In this section we use Theorem \ref{main} and results from algebraic geometry to globally identify the set of Weyl connections on the $2$-sphere whose geodesics are the great circles. 

\subsection{Cartan's connection in the flat case} Consider the \textit{projective} $2$\textit{-sphere} $\mathbb{S}^2=\left(\R^3\setminus\{0\}\right)/\R^{+}$ equipped with its ``outwards'' orientation and projective structure $[\nabla_0]$ whose geodesics are the ``great circles'' $\mathbb{S}^1 \subset\mathbb{S}^2$, i.e.~subspaces of the form $E\cap \mathbb{S}^2$ for some $2$-plane $E\subset \R^3$. Let $(\pi : B\to \mathbb{S}^2,\theta)$ denote Cartan's structure bundle of $(\mathbb{S}^2,[\nabla_0])$. Note that $\mathrm{SL}(3,\R)$ is a right principal $H$-bundle over $\mathbb{S}^2$ with base-point projection 
$$
\tilde{\pi} : \mathrm{SL}(3,\R) \to \mathbb{S}^2, \quad g=(g_0,g_1,g_2) \mapsto [g_0],
$$
where an element $g \in \mathrm{SL}(3,\R)$ is written as three column vectors $(g_0,g_1,g_2)$. Now let $\psi : \mathrm{SL}(3,\R) \to B$ be the map which associates to $g \in \mathrm{SL}(3,R)$ the $2$-frame generated by the map $f_g$ where
$$
f_g : \R^2 \to \mathbb{S}^2,\quad x \mapsto \left[g\cdot \left(\begin{array}{c} x \\ 1 \end{array}\right)\right].$$ 
It turns out that $\psi$ is an $H$-bundle isomorphism which pulls-back $\theta$ to the Maurer Cartan form $\omega$ of $\mathrm{SL}(3,\R)$. In particular since $\d\omega+\omega\wedge\omega=0$, the functions $L_i : \mathrm{SL}(3,\R) \to \R$ vanish. Oriented projective surfaces for which the functions $L_i$ vanish are called \textit{projectively flat} or simply flat. 

\subsection{The bundle of conformal inner products in the flat case}
In the canonical flat case $(\mathbb{S}^2,[\nabla_0])$ the bundle of conformal inner products $\rho : \mathcal{C}(\mathbb{S}^2) \to\mathbb{S}^2$ can be identified explicitly. Let $[b] \in \mathcal{C}(\mathbb{S}^2)$ be a conformal inner product at $[x]=\rho([b])$. Identify $T_{[x]}\mathbb{S}^2$ with $x^{\perp} \subset\mathbb{R}^3$ where $\perp$ denotes the orthogonal complement of $x$ in $\R^3$ with respect to the Euclidean standard metric. Let $(v_1,v_2) \in \R^3\times \R^3$ be a positively oriented conformal basis for $[b]$. The pair $(v_1,v_2)$ is unique up to a transformation of the form  
$$
\left(r\left(v_1\cos\varphi -v_2\sin\varphi\right),r\left(v_1\sin\varphi+v_2\cos\varphi\right)\right)
$$
for some $r \in \R^+$ and $\varphi \in [0,2\pi]$. Thus we may uniquely identify $[b]$ with an element $re^{\i\varphi}\left(v_1+\i v_2\right)$ 
in $\mathbb{CP}^2$. Note that the image of the standard embedding $\iota : \mathbb{RP}^2 \to \mathbb{CP}^2$ precisely consists of those elements $[z] \in \mathbb{CP}^2$ for which $\Re(z)$ and $\Im(z)$ are linearly dependent. It is easy to verify that the just described map is a diffeomorphism $\mathcal{C}(\mathbb{S}^2) \to \cip$ which will be denoted by $\psi$. Thus $\rho_0=\rho \circ \psi^{-1} :\cip \to \mathbb{S}^2$ makes $\cip$ into a $D^2$-bundle over $\mathbb{S}^2$ whose projection map is explicitly given by 
$$
\rho_0 : \cip \to \mathbb{S}^2,\quad [z] \mapsto [\Re(z)\wedge \Im(z)].
$$

\begin{proposition}\label{biholo}
For $(M,[\nabla])=(\mathbb{S}^2,[\nabla_0])$ there exists a biholomorphic fibre bundle isomorphism $\varphi : \mathcal{C}(M) \to \cip$ covering the identity on $\mathbb{S}^2$.   
\end{proposition}
\begin{proof} 
Suppose there exists a smooth surjection $\lambda : \mathrm{SL}(3,\R) \to \cip$ whose fibres are the $C$-orbits and which pulls back the $(1,\!0)$-forms of $\cip$ to linear combinations of $\alpha_1$ and $\alpha_2$. Then it is easy to check that the map $\varphi=\lambda \circ \nu^{-1} : \mathcal{C}(\mathbb{S}^2) \to \cip$ is well defined and has all the desired properties. Consider the smooth map $\tilde{\lambda} : B=\mathrm{SL}(3,\R) \to \C^3$ given by 
$$
\left(g_0 \; g_1 \; g_2\right) \mapsto g_0\wedge \left(g_1+\i g_2\right).
$$
 We have
\begin{equation}\label{xyz}
\tilde{\lambda} \circ R_{h_{a,b}} = \det a^{-1} \left(1,\i\right)\cdot a^t \cdot \left(\begin{array}{c}\Re(\tilde{\lambda}) \\ \Im(\tilde{\lambda})\end{array}\right)
\end{equation}
and 
\begin{equation}\label{xyz2}
\d\tilde{\lambda}=\i g_1 \wedge g_2\left(\theta^1_0+\i\theta^2_0\right)+\Im(\tilde{\lambda})\theta^2_1-\Re(\tilde{\lambda})\theta^2_2+\i \Re(\tilde{\lambda})\theta^1_2-\i\Im(\tilde{\lambda})\theta^1_1. 
\end{equation}
Denote by $q : \C^3\setminus\left\{0\right\} \to \mathbb{CP}^2$ the quotient projection, then \eqref{xyz} implies that $\lambda=q \circ \tilde{\lambda}$ is a smooth surjection onto $\cip$ whose fibres are the $C$-orbits. Moreover it follows with \eqref{xyz2} and straightforward computations that $\lambda$ pulls back the $(1,\! 0)$-forms of $\cip$ to linear combinations of $\alpha_1, \alpha_2$.
\end{proof}

\subsection{Weyl connections on $\mathbb{S}^2$ and smooth quadrics $C \subset \mathbb{CP}^2$} Theorem \ref{main} and Proposition \ref{biholo} now allow to prove: 

\begin{corollary}\label{quadrics}
The Weyl connections on the $2$-sphere whose unparametrised geo\-desics are the great circles are in one-to-one correspondence with the smooth quad\-rics (i.e. smooth algebraic curves of degree $2$) $\mathcal{C} \subset \mathbb{CP}^2$ without real points. 
\end{corollary}
\begin{remark}
The proof can be adapted from~\cite[Theorem 9]{MR1466165}.\footnote{The reason for this is a duality between certain Finsler and Weyl surfaces as reported by Robert L. Bryant in his talk ``Aufwiedersehen Surfaces, revisited'' at the ICM 2006.} The proof given here relies on Theorem \ref{main} and Proposition \ref{biholo}. Another proof could be given by using results from~\cite{MR1979367,lebrunmason}. 
\end{remark} 
\begin{proof}[Proof of Corollary \ref{quadrics}]
Suppose $\nabla$ is a Weyl connection for some conformal structure $[g]$ on $\mathbb{S}^2$ whose geodesics are the great circles. Then by Theorem \ref{main} and Proposition \ref{biholo} $[g] : \mathbb{S}^2 \to \cip$ is a holomorphic curve and hence by Chow's Theorem $\mathcal{C}=[g](\mathbb{S}^2)\subset \cip$ is a smooth algebraic curve whose genus is $0$, since its the image of the $2$-sphere under a section of a fibre bundle. Note that by standard results of algebraic geometry, the genus $g$ and degree~$d$ of a smooth plane algebraic curve satisfy the relation $g=(d-1)(d-2)/2$. It follows that $\mathcal{C}$ is either a line or a quadric. Since every line in $\mathbb{CP}^2$ has a real point, $\mathcal{C}$ must be a quadric. 

Conversely let $\mathcal{C}\subset \mathbb{CP}^2$ be a smooth quadric without real points. In order to show that $\mathcal{C}$ is the image of a smooth section of $\rho : \cip \to \mathbb{S}^2$,  which is a holomorphic curve, it is sufficient to show that $\rho_0\vert_\mathcal{C} : \mathcal{C} \to \mathbb{S}^2$ is a diffeomorphism. The fibre of $\rho : \cip \to \mathbb{S}^2$ at $u=[(u_1,u_2,u_3)^t] \in \mathbb{S}^2$ is an open subset of the real line $u_1z_1+u_2z_2+u_3z_3=0$. Smoothness of $\mathcal{C}$ implies that $\mathcal{C}$ cannot contain that line as a component and since a quadric in $\mathbb{CP}^2$ without real points does not have any real tangent lines, it follows from Bezout's theorem that $\mathcal{C}$ intersects that line transversely in two distinct points. Clearly the line intersects $\mathcal{C}$ in $\mathbb{CP}^1\setminus \mathbb{RP}^1$ which has two connected components: $\rho^{-1}(u)\cup \rho^{-1}(-u)$. The intersection of the quadric with the line consists of one point in each of these components. It follows that $\rho\vert_\mathcal{C} : \mathcal{C} \to \mathbb{S}^2$ is a submersion and hence by the compactness of~$\mathcal{C}$ and~$\mathbb{S}^2$ a covering map which is at most $2$-to-$1$. But since $\mathcal{C}$ is diffeomorphic to $\mathbb{S}^2$,~$\rho$ must restrict to~$\mathcal{C}$ to be a diffeomorphism onto $\mathbb{S}^2$. Hence $(\rho\vert_\mathcal{C})^{-1} : \mathbb{S}^2 \to \cip$ is a smooth section and holomorphic curve whose image is $\mathcal{C}$. According to Theorem \ref{main} this section determines a unique Weyl connection on $\mathbb{S}^2$ whose geodesics are the great circles.
\end{proof}

\subsubsection*{Acknowledgments}
This article is based on the authors doctoral thesis; in this connection the author would like to thank his adviser Norbert Hungerb\"uhler. Also, the author is very grateful to Robert L. Bryant for several stimulating conversations on the topic of this paper. Furthermore, the author would like to thank the referee for his detailed report and valuable suggestions.

\end{document}